\documentclass[a4paper,leqno]{article}
\usepackage{amsmath,amsthm}
\usepackage{mathpazo}
\usepackage[normal]{boilerart}

\newcommand{\Sch}{\categ{Sch}}

\newcommand{\an}{\ensuremath{\textit{an}}}
\newcommand{\coh}{\ensuremath{\textit{coh}}}

\newcommand{\perf}{\ensuremath{\textit{perf}}}

\newcommand{\zar}{\ensuremath{\textit{zar}}}

\DeclareMathOperator{\Gal}{Gal}

\title{On {G}alois categories \& perfectly reduced schemes}

\author{Clark Barwick\thanks{I thank Peter Haine for sharing his many insights about this material. I am also grateful to the Isaac Newton Institute in Cambridge, whose hospitality I enjoyed as I completed this work.}}

\date{}

\begin{document}

\maketitle


Let $X$ be a scheme.\footnote{We only work with coherent schemes, which out of indolence we just call \emph{schemes}.} A point $x=\Spec\kappa(x)\to X$ of $X$ has an image Zariski point $x_0\in X^{\zar}$ with residue field $\kappa(x_0)\subseteq\kappa(x)$. Let us say that $x$ is a \emph{geometric} point when $\kappa(x)$ is a separable closure of $\kappa(x_0)$. Geometric points constitute a category \cite[Exposé VIII, \S 7]{MR50:7131}, which we call the \emph{Galois category} $\Gal(X)$. The morphisms $x\to y$ are \emph{specialisations} $x\leftsquigarrow y$ -- i.e., natural transformations between the corresponding morphisms of topoi $x_{\ast}\ot y_{\ast}$ (or, if you prefer, $x^{\ast}\to y^{\ast}$). In other words, $\Gal(X)$ is the category of points of the étale topos of $X$. It is a $1$-category in which every endomorphism is an automorphism. It comes equipped with a profinite topology; that is, $\Gal(X)$ is a category object in Stone topological spaces.

The Galois category also comes equipped with a conservative functor $\Gal(X)\to X^{\zar}$, whose target is a poset under specialisation; this functor is continuous for the profinite topologies.\footnote{The topological space $X^{\zar}$ is a \emph{spectral topological space}, which is the same thing as a profinite poset.} Accordingly, $X^{\zar}$ is the poset of isomorphism classes of objects of $\Gal(X)$.

The profinite category $\Gal(X)$ is determined by the étale topos of $X$, but it also determines it; in fact, if you're a hyperpolyglot, you can probably deduce this already from Makkai's Strong Conceptual Completeness Theorem \cite{MR900266}. We took \cite{exodromy} an explicit approach that showed that étale sheaves on $X$ `are' continuous representations $\Gal(X)$, generalising the usual equivalence between the étale cohomology and the Galois cohomology of a field.

If $X^{\zar}\to P$ is a finite constructible stratification of a scheme, then the \emph{Galois $\infty$-category} $\Gal(X/P)$ is what you get by localising (in the wholesome $\infty$-categorical sense) the specialisations that occur within any single stratum. The result is a profinite $\infty$-category with a conservative functor to $P$ -- what we called a \emph{profinite $P$-stratified space} \cite{exodromy}. It is the $\infty$-category of points of the \emph{$P$-stratified $\infty$-topos}. In the extreme case, when $P$ is the trivial poset, $\Gal(X/\ast)$ is the profinite étale homotopy type. Hence $\Gal(X)$ is a complete delocalisation of the étale homotopy type.

When you view $\Gal(X)$ through this lens, you get to interpret it as a profinite stratified space whose underlying space is the profinite étale homotopy type $\Gal(X/\ast)$. Each irreducible closed subscheme $Z\subseteq X$ identifies the closure $[Z]$ of a stratum within $X$. If $Z\subseteq W$ are two irreducible closed subschemes of $X$, then the space of sections of $\Gal(X)\to X^{\zar}$ over the edge $\eta_Z\to\eta_W$ of the generic points is the deleted tubular neighbourhood\footnote{This is literally the étale homotopy type of the \emph{oriented fibre product} $\eta_Z\mathbin{\overleftarrow{\times}}_X\eta_W$.} of $[Z]$ in $[W]$. This stratified space is a stratified $1$-type: the strata and deleted tubular neighbourhoods are all $K(\pi,1)$'s.

\begin{exm*}[Fields] If $k$ is a field, then a choice of a separable closure of $k$ identifies an equivalence $\Gal(\Spec k)\simeq BG_k$, where $G_k$ is the absolute Galois group of $k$.
\end{exm*}

\begin{exm*}[Knots and primes] If $A$ is a number ring with fraction field $K$, then $\Gal(\Spec A)$ is a category with (isomorphism classes of) objects the prime ideals of $A$. For each nonzero prime ideal $\mathfrak{p}\in\Spec A$, the automorphisms of $\mathfrak{p}$ can be identified with the absolute Galois group $G_{\kappa(\mathfrak{p})}$ of the finite field $\kappa(\mathfrak{p})$. Thus the étale homotopy type of $\Spec A$ is stratified by the various closed strata, each of which is an embedded circle -- i.e., a knot $BG_{\kappa(\mathfrak{p})}$. The open complement of each $BG_{\kappa(\mathfrak{p})}$ is a $BG_{\mathfrak{p}}$, where $G_{\mathfrak{p}}\coloneq\pi_1(\Spec A\smallsetminus\mathfrak{p})$ is the automorphism group of the maximal Galois extension of $K$ that is ramified at most only at $\mathfrak{p}$ and the infinite primes. Enveloping each knot is a tubular neighbourhood, given by $\Gal(\Spec A_{\mathfrak{p}}^{\textit{sh}})$ (sh=strict henselisation), so that the deleted tubular neighbourhood of $BG_{\kappa(\mathfrak{p})}$ is a $BG_{K_{\mathfrak{p}}}$.
\end{exm*}

\begin{exm*}[Analytification] If $X$ is a finite type $F$-scheme, where $F$ is $\CC$, $\RR$, or any nonarchimedean field, then there is an associated $X^{\textit{an}}$ analytic space, which admits a profinite stratification by $X^{\zar}$. The category $\Gal(X)$ is the profinite completion of the exit-path $\infty$-category of $X^{\an}$ with this stratification. (We proved this over $\CC$ \cite[Proposition 13.15 \& Corollary 13.16]{exodromy}, but the same proof will work any time you have access to an Artin Comparison Theorem, which you do in these situations; see \cite{MR1357745}.)
\end{exm*}

The \emph{perfectly reduced schemes} of the title are schemes taken up to universal homeomorphism (\Cref{dfn:perfectlyreduced}). Grothendieck's \emph{invariance topologique} of the étale topos \cite[Exposé VIII, 1.1]{MR50:7131} ensures that the \emph{only} kinds of schemes Galois categories can hope to capture in their entirety are the perfectly reduced schemes. This note is the suggestion of a recognition principle that flows in the opposite direction; that is, we aim to read off facts about perfectly reduced schemes from their Galois categories. Our goal is a dictionary between the geometric features of a perfectly reduced scheme (or morphism of such) and the categorical properties of its Galois category (or functor of such); the gnomic section titles are the first few entries in this dictionary.

\tableofcontents


\section{Open = cosieve \& closed = sieve} Let us begin with the obvious.

\begin{prp} A monomorphism $U\inclusion X$ of schemes is an \emph{open immersion} if and only if the induced functor $\Gal(U)\to\Gal(X)$ is equivalent to the inclusion of a cosieve.

Dually, a monomorphism $Z\inclusion X$ of schemes is a \emph{closed immersion} if and only if $\Gal(Z)\to\Gal(X)$ is equivalent to the inclusion of a sieve.
\end{prp}

An \emph{interval} in an $\infty$-category $C$ is a full subcategory $D\subseteq C$ such that a morphism $P\to Q$ of $D$ factors through an object $R$ of $C$ only if $R$ lies in $D$.
\begin{cor} A monomorphism $W\inclusion X$ of schemes is a \emph{locally closed immersion} if and only if the induced functor $\Gal(W)\to\Gal(X)$ is equivalent to the inclusion of an interval.
\end{cor}

\begin{cor} A scheme $X$ is local if and only if $\Gal(X)$ contains a \emph{weakly initial} object -- i.e., an object from which every object receives a morphism. Dually, a scheme $X$ is irreducible if and only if $\Gal(X)$ contains a \emph{weakly terminal} object -- i.e., an object to which every object sends a morphism.
\end{cor}

\begin{nul} For any scheme $X$ and any point $x_0\in X^{\zar}$, the Galois category of the localisation is the fibre product
\[
\Gal(X_{(x_0)})\simeq\Gal(X)\times_{X^{\zar}}X^{\zar}_{x/} \rlap{\ .}
\]
Dually, for any point $y_0\in X^{\zar}$, the Galois category of the closure $X^{(y_0)}$ of $y_0$ (with the reduced subscheme structure, say) is the fibre product
\[
\Gal(X^{(y_0)})\simeq\Gal(X)\times_{X^{\zar}}X^{\zar}_{/y} \rlap{\ .}
\]
\end{nul}


\section{Strict localisation = undercategory \& strict normalisation = overcategory} 

\begin{ntn} If $x\to X$ is a point of a scheme $X$, then we write $O^h_{X,x_0}$ for the henselisation of the local ring $O_{X,x_0}$, and we write $O^h_{X,x}\supseteq O^h_{X,x_0}$ for the unique extension of henselian local rings that on residue fields reduces to the field extension $\kappa\supseteq\kappa(x_0)$, where $\kappa$ is the separable closure of $\kappa(x_0)$ in $\kappa(x)$. We will also write
\[
X_{(x)}\coloneq\Spec O^h_{X,x} \rlap{\ .}
\]
We call $X_{(x)}$ the \emph{localisation of $X$ at $x$}. It is the limit of the factorisations $x\to U\to X$ in which $U\to X$ is étale.

If $x\to X$ is a geometric point, then $O^h_{X,x}$ is the strict henselisation of $O_{X,x_0}$, and $X_{(x)}$ is the strict localisation of $X$ at $x$.

Dually, if $y\to X$ is a point, then we write $X^{(y_0)}$ for the reduced subscheme structure on the Zariski closure of $y_0$, and we write $X^{(y)}$ for the normalisation of $X^{(y_0)}$ under $\Spec\kappa$, where $\kappa$ is the separable closure of $\kappa(y_0)$ in $\kappa(y)$. We call $X^{(y)}$ the \emph{normalisation of $X$ at $y$}.

If $y\to X$ is a geometric point, then we call $X^{(y)}$ the \emph{strict normalisation of $X$ at $y$}. It is the limit of the factorisations $y\to Z\to X$ in which $Z\to X$ is finite.
\end{ntn}

\begin{nul} Stefan Schröer \cite{MR3649361} has brought us \emph{totally separably closed} schemes, which are integral normal schemes whose function field is separably closed. In other words, a totally separably closed scheme is one of the form $X^{(y)}$ for some geometric point $y\to X$. (In the language of Schröer, $X^{(y)}$ is the total separable closure of the Zariski closure of $y_0$ -- with the reduced subscheme structure -- under $y$.) Schröer has shown that this class of schemes has a number of curious properties:
\begin{itemize}
\item If $Z$ is totally separably closed, then for any point $z_0\in Z^{\zar}$, the local ring $O_{Z,z_0}$ is strictly henselian \cite[Proposition 2.6]{MR3649361}.
\item If $Z$ is totally separably closed, then the étale topos and the Zariski topos of $Z$ coincide, so that $\Gal(Z)\simeq Z^{\zar}$ \cite[Corollary 2.5]{MR3649361}. In other words, $\Gal(Z)$ is a profinite poset with a terminal object.
\item If $Z$ is totally separably closed and $W$ is irreducible, then any integral morphism $W\to Z$ is radicial \cite[Lemma 2.3]{MR3649361}. Thus any integral surjection $W\to Z$ is a universal homeomorphism.
\item If $Z$ is totally separably closed, then the poset $\Gal(Z)\simeq Z^{\zar}$ has all finite nonempty joins \cite[Theorem 2.1]{tsccontractions}.
\end{itemize}
\end{nul}

Here now is the basic observation, which follows more or less immediately from the limit descriptions of the strict localisation and the strict normalisation:
\begin{prp} Let $X$ be a scheme, and let $x\to X$ and $y\to X$ be two geometric points thereof. The following profinite sets are in (canonical) bijection:
\begin{itemize}
\item the set $\Map_{\Gal(X)}(x,y)$ of morphisms $x\to y$ in $\Gal(X)$;
\item the set $\Mor_X(y,X_{(x)})$ of lifts of $y$ to the strict localisation $X_{(x)}$;
\item the set $\Mor_X(x,X^{(y)})$ of lifts of $y$ to the strict normalisation $X^{(y)}$.
\end{itemize}
\end{prp}
\noindent We may thus describe the over- and undercategories of Galois categories:
\begin{cor} Let $X$ be a  scheme, and let $x\to X$ and $y\to X$ be two geometric points thereof. Then we have
\[
\Gal(X)_{x/}\simeq\Gal(X_{(x)})\text{\quad and\quad}\Gal(X)_{/y}\simeq\Gal(X^{(y)}) \rlap{\ .}
\]
\end{cor}
\noindent See also \cite[Exposé VIII, Corollaire 7.6]{MR50:7131}, where the first sentence is proved.
\begin{cor} Let $X$ be a scheme. Then $\Gal(X)$ is equivalent to both of the following full subcategories of $X$-schemes:
\begin{itemize}
\item the one spanned by the strict localisations of $X$, and
\item the one spanned by the strict normalisations of $X$.
\end{itemize}
\end{cor}
\noindent Since $\Gal(X^{(y)})\simeq X^{(y),\zar}$, it follows that Galois categories are of a very particular sort:
\begin{cor} Let $X$ be a  scheme. For any geometric point $y\to X$, the overcategory $\Gal(X)_{/y}$ is a profinite poset with all finite nonempty joins. In particular, every morphism of $\Gal(X)$ is a monomorphism. 
\end{cor}

\begin{dfn} Let $X$ be a scheme. Then a \emph{witness} is a totally separably closed valuation ring $V$ and a morphism $\gamma\colon\Spec V\to X$. If $p_0$ is the initial object of $\Gal(V)$ and $p_{\infty}$ is the terminal object of $\Gal(V)$, then we say that $\gamma$ \emph{witnesses} the map $\gamma(p_0)\to\gamma(p_{\infty})$ of $\Gal(X)$.
\end{dfn}

\begin{nul} Any morphism $x\to y$ of $\Gal(X)$ has a witness: you can always find a local morphism $\Spec V\to (X^{(y)})_{(x)}$ that induces an isomorphism of function fields.
\end{nul}


\section{Universal homeomorphism = equivalence}

Now we arrive at a sensitive question: under which circumstances does a morphism of  schemes induce an equivalence of étale topoi or, equivalently, of Galois categories? The well-known theorem here is Grothendieck's \emph{invariance topologique} of the étale topos \cite[Exposé VIII, 1.1]{MR50:7131}, which states that a universal homeomorphism induces an equivalence on étale topoi. Let us reprove this result with the aid of Galois categories; this will also provide us with a partial converse.

\begin{prp}\label{prp:radicialintermsofGal} Let $f\colon X\to Y$ be a morphism of  schemes. If $f$ is radicial, then every fibre of $\Gal(X)\to\Gal(Y)$ is either empty or a singleton.\footnote{By \emph{singleton} we mean \emph{contractible groupoid}.} Conversely, if $f$ is of finite type, and if every fibre of $\Gal(X)\to\Gal(Y)$ is either empty or a singleton, then $f$ is radicial.
\end{prp}

\begin{proof} If $f$ is radicial, then the map $X^{\zar}\to Y^{\zar}$ is an injection, and for any point $x_0\in X^{\zar}$, the map $BG_{\kappa(x_0)}\to BG_{\kappa(f(x_0))}$ on fibres is an equivalence since $\kappa(f(x_0))\subseteq\kappa(x_0)$ is purely inseparable. So for any geometric point $y$ with image $y_0$, the fibre over $y$ is a singleton.

Conversely, if $f$ is of finite type, and if every fibre of $\Gal(X)\to\Gal(Y)$ is either empty or a singleton, then certainly the map $X^{\zar}\to Y^{\zar}$ is an injection, whence $f$ is in particular quasifinite. For any point $x_0\in X^{\zar}$, the fibres of the map $BG_{\kappa(x_0)}\to BG_{\kappa(f(x_0))}$ are each a singleton, whence it is an equivalence. Now since $\kappa(f(x_0))\subseteq\kappa(x_0)$ is a finite extension, it is purely inseparable.
\end{proof}

\begin{exm} The finite type hypothesis in the second half of \Cref{prp:radicialintermsofGal} is of course necessary, as any nontrivial extension $E\subset F$ of separably closed fields induces the identity on trivial Galois categories.
\end{exm}

\begin{cor}\label{cor:radicialsurjintermsofGal} Let $f\colon X\to Y$ be a morphism of  schemes. If $f$ is radicial and surjective, then every fibre of $\Gal(X)\to\Gal(Y)$ is a singleton. Conversely, if $f$ is of finite type, and if every fibre of $\Gal(X)\to\Gal(Y)$ is a singleton, then $f$ is radicial and surjective.
\end{cor}

The following is the Valuative Criterion, along with a simple argument \stacks{03K8} that allows one to extend the fraction field of the valuation ring therein.
\begin{lem}\label{lem:witnessesofunivclosed} Let $f\colon X\to Y$ be a morphism of  schemes. Then the following are equivalent.
\begin{itemize}
\item The morphism $f$ is universally closed.
\item For any witness $\gamma\colon\Spec V\to Y$ and any diagram
\[
\begin{tikzcd}
	\Spec K \arrow[r] \arrow[d] & X \arrow[d, "f" right] \\ 
	\Spec V \arrow[r, "\gamma" below] & Y
\end{tikzcd}
\]
in which $K$ is the fraction field of $V$, there exists a lift $\overline{\gamma}\colon\Spec V\to X$.
\end{itemize}
\end{lem}

\begin{rec} A functor $f\colon C\to D$ is said to be a \emph{right fibration} if and only if, for any object $x\in C$, the induced functor $C_{/x}\to D_{/f(x)}$ is an equivalence of categories. In this case, one may say that $C$ is a \emph{category fibred in groupoids over $D$}. For any such right fibration, there is a diagram $F$ of groupoids indexed on $D^{\op}$ such that $C$ is the Grothendieck construction of $F$.

Dually, $f$ is a \emph{left fibration} if and only if $f^{\op}$ is a right fibration, so that for any object $x\in C$, the induced functor $C_{x/}\to D_{f(x)/}$ is an equivalence of categories.
\end{rec}

\begin{prp}\label{prp:integralintermsofGal} Let $f\colon X\to Y$ be a morphism of schemes. If $f$ is an integral morphism, then $\Gal(X)\to\Gal(Y)$ is a right fibration. Conversely, if $\Gal(X)\to\Gal(Y)$ is a right fibration, then $f$ is universally closed.
\end{prp}

\begin{proof} Assume that $f$ is integral. Then for every geometric point $x\to X$, the induced morphism $X^{(x)}\to Y^{(f(x))}$ is also integral, and by \cite[Lemma 2.3]{MR3649361}, it is radicial as well. Hence at the level of Zariski topological spaces, $X^{(x),\zar}\to Y^{(f(x)),\zar}$ is an inclusion of a closed subset; since source and target are each irreducible, and the inclusion carries the generic point to the generic point, it is a homeomorphism. (In fact, $X^{(x)}\to Y^{(f(x))}$ is a universal homeomorphism.) Thus
\[
\Gal(X)_{/x}\simeq\Gal(X^{(x)})\simeq X^{(x),\zar}\to Y^{(f(x)),\zar}\simeq\Gal(Y^{(f(x))})\simeq\Gal(Y)_{/f(x)}
\]
is an equivalence, whence $\Gal(X)\to\Gal(Y)$ is a right fibration.

Conversely, assume that $f$ is of finite type and that $\Gal(X)\to\Gal(Y)$ is a right fibration. We employ \Cref{lem:witnessesofunivclosed} to show that $f$ is universally closed; consider a witness $\gamma\colon\Spec V\to Y$ along with a diagram
\[
\begin{tikzcd}
	\Spec K \arrow[r, "\xi" above] \arrow[d] & X \arrow[d, "f" right] \\ 
	\Spec V \arrow[r, "\gamma" below] & Y
\end{tikzcd}
\]
in which $K$ is the fraction field of $V$. Let $\psi\colon y\to f(\xi)$ be the morphism of $\Gal(Y)$ witnessed by $\gamma$, and let $\phi\colon x\to\xi$ be a lift thereof to $\Gal(X)$. We obtain a square
\[
\begin{tikzcd}
	O_{Y,y}^{\textit{sh}} \arrow[r, "\gamma" above] \arrow[d] & V \arrow[d] \\ 
	O_{X,x}^{\textit{sh}} \arrow[r, "\xi" below] & K \rlap{\ ,}
\end{tikzcd}
\]
and since $O_{Y,y}^{\textit{sh}}\to O_{X,x}^{\textit{sh}}$ is local, we obtain a lift $\overline{\gamma}\colon O_{X,x}^{\textit{sh}}\to V$, as required.
\end{proof}

A universal homeomorphism is a morphism that is radicial, surjective, and universally closed. An equivalence of categories is a right fibration with fibres contractible groupoids. We thus deduce:
\begin{prp} Let $f\colon X\to Y$ be a morphism of schemes. If $f$ is a universal homeomorphism, then $\Gal(X)\to\Gal(Y)$ is an equivalence. Conversely, if $f$ is of finite type, and if $\Gal(X)\to\Gal(Y)$ is an equivalence, then $f$ is a universal homeomorphism (which is necessarily finite).
\end{prp}

\section{Interlude: perfectly reduced schemes}\label{sec:interlude}

A reduced scheme receives no nontrivial nilimmersions; a \emph{perfectly reduced} scheme receives no nontrivial universal homeomorphisms. This is in fact a local condition that can be expressed in very concrete terms:
\begin{prp}\label{prp:trig} The following are equivalent for a scheme $X$.
\begin{itemize}
\item There exists an affine open covering $\{\Spec A_i\}_{i\in I}$ of $X$ such that for every $i\in I$, the following conditions obtain:
\begin{itemize}
\item for any $f,g\in A_i$, if $f^2=g^3$, then there is a unique $h\in A_i$ such that $f=h^3$ and $g=h^2$; and
\item for any prime number $p$ and any $f,g\in A_i$, if $f^p=p^pg$, then there is a unique element $h\in A_i$ such that $f=ph$ and $g=h^p$.
\end{itemize} 
\item If $X'$ is a reduced scheme and $f\colon X'\to X$ is a universal homeomorphism, then $f$ is an isomorphism.
\end{itemize}
\end{prp}

\begin{dfn}\label{dfn:perfectlyreduced} A  scheme that enjoys one and therefore both of the conditions of \Cref{prp:trig} is said to be \emph{perfectly reduced} or -- in the parlance of \cite[Appendix B]{MR2679038} and \stacks{0EUL} -- \emph{absolutely weakly normal}.

Let us write $\Sch_{\perf}\subset\Sch_{\coh}$ for the full subcategory of schemes spanned by the perfectly reduced schemes.
\end{dfn}

\begin{nul} To express this differently, let us define a family of reference universal homeomorphisms. First, let $\Upsilon$ denote the cuspidal cubic
\[
\Upsilon\coloneq\Spec\ZZ[u,v]/(u^2-v^3) \rlap{\ .}
\]
The normalisation $\rho\colon\AA^1_{\ZZ}\to\Upsilon$ defined by the equations $u=t^3$ and $v=t^2$ is a universal homeomorphism. Next, for any prime number $p$, set
\[
Z_p\coloneq\Spec\ZZ[y,z]/(y^p-p^pz) \rlap{\ .}
\]
The normalisation $\tau_p\colon\AA^1_{\ZZ}\to Z_p$ defined by the equations $y=px$ and $z=x^p$ is a universal homeomorphism. \Cref{prp:trig} states that a  scheme $X$ is perfectly reduced if and only if every point $x\in X$ is contained in a Zariski open neighbourhood $U\subseteq X$ such that the map
\[
\Mor(U,\AA^1_{\ZZ})\to\Mor(U,\Upsilon)
\]
is a bijection, and for any prime number $p$, the map
\[
\Mor(U,\AA^1_{\ZZ})\to\Mor(U,Z_p)
\]
is a bijection.
\end{nul}

\begin{nul} Any (quasicompact) open subscheme of a perfectly reduced scheme is perfectly reduced. A reduced $\QQ$-scheme is perfectly reduced if and only if it is \emph{seminormal}. A reduced $\FF_{\!p}$-scheme is perfectly reduced if and only if the Frobenius morphism is an isomorphism.
\end{nul}

\begin{prp}[\protect{\cite[Proposition 14.5]{exodromy}}] The inclusion $\Sch_{\perf}\inclusion\Sch_{\coh}$ admits a right adjoint $\goesto{X}{X_{\perf}}$, which exhibits $\Sch_{\perf}$ as the colocalisation of $\Sch_{\coh}$ along the class of universal homeomorphisms. In particular, the counit $X_{\perf}\to X$ is the initial object in the category of universal homeomorphisms to $X$. We call $X_{\perf}$ the \emph{perfection} of $X$.
\end{prp}

\begin{nul} For reduced $\QQ$-schemes, the perfection is the seminormalisation \stacks{0EUT}. For reduced $\FF_{\!p}$-schemes $X$ the perfection is the limit of $X$ over powers of the Frobenius, as usual.
\end{nul}

\begin{dfn} A \emph{topological morphism} from a scheme $X$ to a scheme $Y$ is an morphism $\phi\colon X_{\perf}\to Y$. If $\phi$ induces an isomorphism $X_{\perf}\equivalence Y_{\!\perf}$, then it is said to be a \emph{topological equivalence} from $X$ to $Y$.
\end{dfn}

\begin{nul} Let $X$ and $Y$ be schemes. Consider the following category $T(X,Y)$. The objects are diagrams 
\[
X\ot X'\to Y
\]
in which $X\ot X'$ is a universal homeomorphism. A morphism
\[
\text{from\qquad}X\ot X'\to Y\text{\qquad to\qquad}X\ot X''\to Y
\]
is a commutative diagram
\[
\begin{tikzcd}[row sep=tiny,column sep=tiny]
& X' \arrow[dl] \arrow[dd] \arrow[dr] & \\ 
X && Y\\
& X'' \arrow[ul] \arrow[ur] &
\end{tikzcd}
\]
in which the vertical morphism is (of necessity) a universal homeomorphism. The nerve of the category $T(X,Y)$ is equivalent to the set $\Mor(X_{\perf},Y)\cong\Mor(X_{\perf},Y_{\!\perf})$ of topological morphisms from $X$ to $Y$.
\end{nul}

\begin{nul} The point now is that $\Gal$, viewed as a functor from $\Sch_{\perf}$ to categories, is \emph{conservative}.
\end{nul}

\begin{dfn} Let $P$ be a property of morphisms of schemes that is stable under base change and composition. We will say that a morphism $f\colon X\to Y$ is \emph{topologically $P$} if and only if it is topologically equivalent to a morphism of schemes $f'\colon X'\to Y'$ with property $P$.
\end{dfn}

\begin{nul} Let $P$ be a property of morphisms of schemes that is stable under base change and composition. The class of topologically $P$ morphisms is the smallest class of morphisms $P^t$ that contains $P$ and satisfies the following condition: for any commutative diagram
\[
\begin{tikzcd}
X \arrow[r, "f" above] \arrow[d, "\phi" left] & Y \arrow[d, "\psi" right]\\ 
X' \arrow[r, "f'" below] & Y'
\end{tikzcd}
\]
in which $\phi$ and $\psi$ are universal homeomorphisms, the morphism $f$ lies in $P^t$ if and only if $f'$ does.

A morphism $f\colon X\to Y$ of perfectly reduced schemes is topologically $P$ precisely when it factors as a universal homeomorphism $X\to X'$ followed by a morphism $X'\to Y$ with property $P$.
\end{nul}

\begin{exm} A morphism $f\colon X\to Y$ of perfectly reduced schemes is topologically radicial, surjective, universally closed, or integral if and only if it is radicial, surjective, universally closed, or integral (respectively).
\end{exm}

\begin{exm} A morphism $f\colon X\to Y$ of perfectly reduced schemes is topologically étale if and only if it is étale. Indeed, if $f'\colon X'\to Y$ is étale, then $X'$ is perfectly reduced \cite[B.6(ii)]{MR2679038}.
\end{exm}


\section{Finite = right fibration with finite fibres} We've already seen that an integral morphism of schemes induces a right fibration of Galois categories and that a morphism that induces a right fibration of Galois categories must be universally closed. Let us complete this picture.

Let us begin with an obvious characterisation of quasifinite morphisms. We will say that a functor has \emph{finite fibres} if each of its fibres is a finite set\footnote{which for our purposes means a finite disjoint union of contractible groupoids}.
\begin{lem} Let $f\colon X\to Y$ be a morphism that is of finite type. Then $f$ is quasifinite if and only if $\Gal(X)\to\Gal(Y)$ has finite fibres.
\end{lem}

Since proper quasifinite morphisms are finite, \Cref{prp:integralintermsofGal} now yields:
\begin{prp}\label{prp:finiteintermsofGal} Let $f\colon X\to Y$ be a morphism that is separated and of finite type. Then $f$ is finite if and only if $\Gal(X)\to\Gal(Y)$ is a right fibration with finite fibres.
\end{prp}


\section{Étale = left fibration with finite fibres}

\begin{prp}\label{prp:etaleintermsofGal} Let $f\colon X\to Y$ be a morphism of schemes. If $f$ is weakly étale, then $\Gal(X)\to\Gal(Y)$ is equivalent to a left fibration. Conversely, if $X$ and $Y$ are perfectly reduced, if $f$ is of finite presentation, and if $\Gal(X)\to\Gal(Y)$ is a left fibration with finite fibres, then $f$ is étale.
\end{prp}

\begin{proof} Assume that $f$ is weakly étale. Then for any geometric point $x\to X$, the morphism $X_{(x)}\to Y_{(f(x))}$ is an isomorphism, whence the functor
\[
\Gal(X)_{x/}\simeq\Gal(X_{(x)})\to\Gal(Y_{(f(x))})\simeq\Gal(Y)_{f(x)/}
\]
is an equivalence, whence $\Gal(X)\to\Gal(Y)$ is a left fibration.

Conversely, assume that $X$ and $Y$ are perfectly reduced, that $f$ is of finite presentation, and that $\Gal(X)\to\Gal(Y)$ is a left fibration with finite fibres. So the functor $\Gal(X)\to\Gal(Y)$ is classified by a continuous functor $\Gal(Y)\to\Set^{\fin}$, which in turn corresponds to a constructible étale sheaf of finite sets on $Y$, which in particular coincides with the sheaf represented by $X$. Since the sheaf represented by $X$ is constructible, there exists an étale map $U\to Y$ and an effective epimorphism $U\to X$ of étale sheaves on $Y$. By descent, $X\to Y$ is étale.
\end{proof}


\section{Finite étale = Kan fibration with finite fibres}

We may as well combine the last two entries in our dictionary.
\begin{rec} A \emph{Kan fibration} is a functor that induces a Kan fibration on nerves. Equivalently, it is a functor that is both a left and right fibration. Equivalently, it is a functor $C\to D$ that is equivalent to the Grothendieck construction applied to a diagram of groupoids indexed on $D^{\op}$ that carries every morphism to an equivalence of groupoids.
\end{rec}

\begin{prp} Let $f\colon X\to Y$ be a morphism of perfectly reduced schemes that is separated and of finite presentation. Then $f$ is finite étale if and only if $\Gal(X)\to\Gal(Y)$ is a Kan fibration with finite fibres.
\end{prp}


\DeclareFieldFormat{labelnumberwidth}{#1}
\printbibliography[keyword=alph]
\addcontentsline{toc}{section}{References} 
\DeclareFieldFormat{labelnumberwidth}{{#1\adddot\midsentence}}
\printbibliography[heading=none, notkeyword=alph]

\end{document}